\documentclass[letterpaper, 10 pt, conference]{ieeeconf}  

\IEEEoverridecommandlockouts                              
\usepackage{amsmath,amssymb,amsfonts}
\usepackage{balance}
\usepackage{color}
\usepackage{mathtools}
\usepackage{pdfrender}

\usepackage[acronym]{glossaries}
\usepackage[pdftex]{graphicx}
\usepackage{subfig}
\usepackage{enumerate}
\usepackage{tabularx}
\usepackage{array,booktabs}
\usepackage{multirow}
\usepackage[linesnumbered,ruled,vlined]{algorithm2e}

\graphicspath{{Pictures/}}

\newtheorem{definition}{Definition}
\newtheorem{proposition}{Proposition}

\makeglossaries
\newacronym{MIQP}{MIQP}{Mixed-Integer Quadratic Programming}
\newacronym{AD}{AD}{Automated Driving}
\newacronym{MILP}{MILP}{Mixed-Integer Linear Programming}
\newacronym{MLD}{MLD}{Mixed-Logical-Dynamical}
\newacronym{MPC}{MPC}{Model Predictive Control}
\newacronym{GS}{GS}{Gauss-Southwell}
\newacronym{GNEP}{GNEP}{Generalized Nash Equilibrium Problem}
\newacronym{MINE}{$\varepsilon$-MINE}{$\varepsilon$-Mixed-Integer Nash Equilibrium}

\newcommand{\R}{\mathbb{R}}

\newcommand{\bs}[1]{\boldsymbol{#1}}
\newcommand{\mc}[1]{\mathcal{#1}}

\newcommand{\figname}{Fig.}

\title{\LARGE \bf
A Mixed-Logical-Dynamical model for Automated Driving on highways
}

\author{
Filippo Fabiani and Sergio Grammatico
\thanks{F. Fabiani is with the Department of Information Engineering, University of Pisa, Italy.
S. Grammatico is with the Delft Center for Systems and Control (DCSC), TU Delft, The Netherlands. E-mail addresses: \texttt{filippo.fabiani@ing.unipi.it}, \texttt{s.grammatico@tudelft.nl}. 
This work was partially supported by NWO under research projects OMEGA (grant n. 613.001.702) and P2P-TALES (grant n. 647.003.003), and by DCSC/3mE/TU Delft under research project Intelligent Autonomous Vehicles.
		\smallskip \newline
	}
}

\begin{document}

\maketitle
\thispagestyle{empty}
\pagestyle{empty}

\begin{abstract}
We propose a hybrid decision-making framework for safe and efficient autonomous driving of selfish vehicles on highways. 
Specifically, we model the dynamics of each vehicle as a Mixed-Logical-Dynamical system and propose simple driving rules to prevent potential sources of conflict among neighboring vehicles. We formalize the coordination problem as a generalized mixed-integer potential game, where an equilibrium solution generates a sequence of mixed-integer decisions for the vehicles that trade off individual optimality and overall safety.
\end{abstract}

\section{Introduction}

\gls{AD} is currently foreseen as the future of road traffic to enhance safety and efficiency. 
Within the system-and-control community, multi-vehicle coordination, motion planning and control for AD has attracted a strong research attention, since it poses relevant engineering challenges,  spanning from fundamental to computational and practical challenges. Providing each vehicle with a high degree of decision-making autonomy is in fact key towards automated road traffic.
From an optimal-control perspective, to autonomously drive vehicles within a complex dynamic environment, several algorithms propose a \gls{MPC} approach \cite{falcone2007predictive}, \cite{glaser2010maneuver}, \cite{kim2014mpc}, or some variants, such as scenario-based \gls{MPC} \cite{cesari2017scenario}, spatial-based \gls{MPC} \cite{plessen2018spatial}, and distributed \gls{MPC} \cite{keviczky2008decentralized, mohseni2017distributed}, as well as multi-layer decision-making frameworks \cite{gray2012predictive}, \cite{noh2018decision}.

The quintessential feature in multi-vehicle driving scenarios is that drivers are selfish decision makers, or agents, that pursue their own individual interests, e.g. minimum travel time or minimum fuel consumption, while sharing the road space-time. To handle the presence of selfish vehicles, the principles of \textit{game theory} have been adopted, first in high-level traffic control \cite{baskar2012traffic}, and more recently in multi-vehicle motion planning \cite{wang2015game, bahram2016game, li2017game}.

\begin{figure}[!t]
	\centering
	\includegraphics[width=0.95\columnwidth]{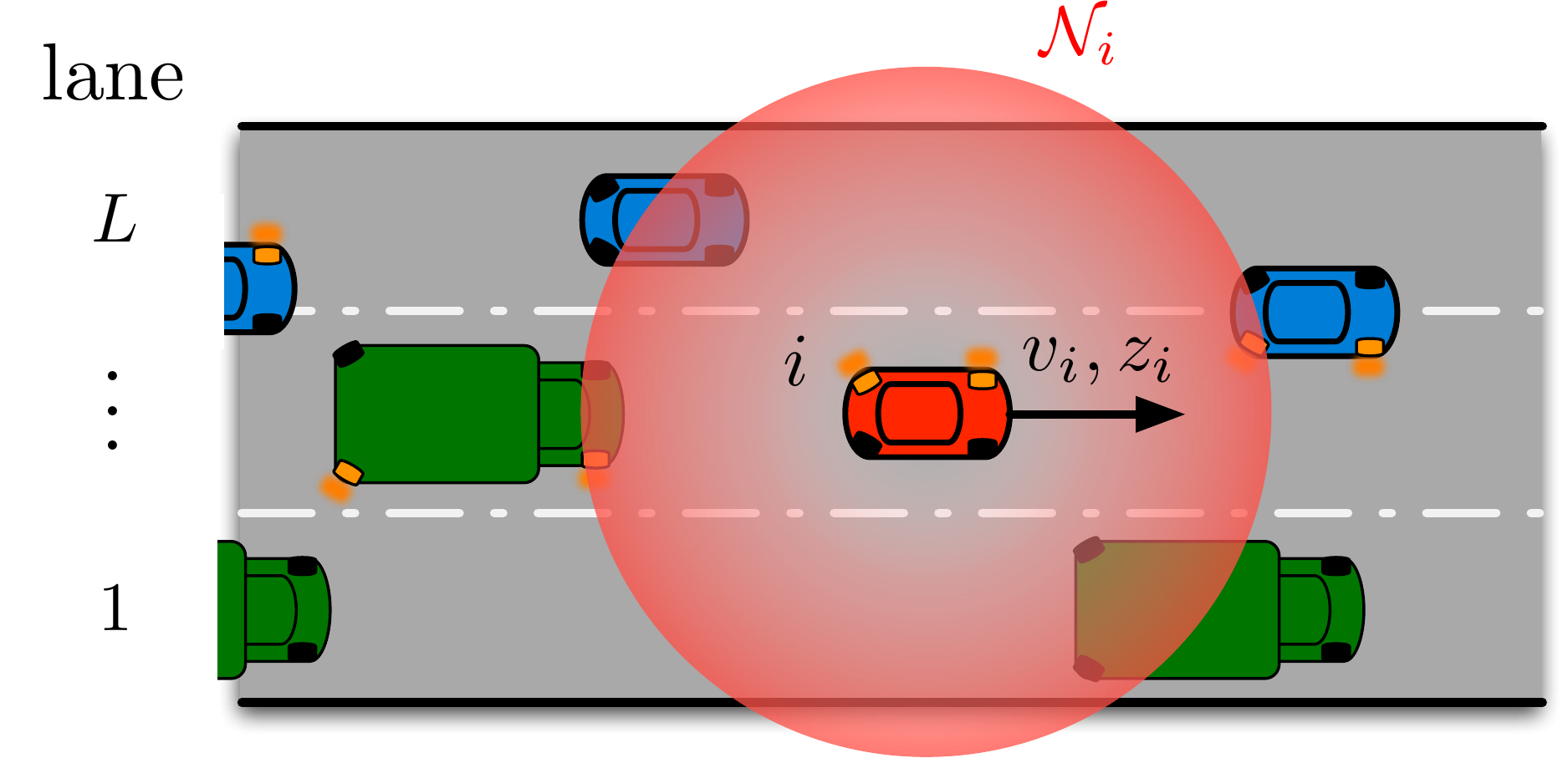}
	\caption{A set of vehicles driving along a highway. 
			}
	\label{fig:AutPlayers}
\end{figure}

In this paper, compared with the referred literature, we consider a general driving scenario on multi-lane highways with multiple vehicles, each with a cost function to be minimized given the driving decisions of the other vehicles, individual constraints, e.g. speed and acceleration limits, and safety-distance constraints. Furthermore, for each vehicle, we embed both continuous and discrete decisions over a prediction horizon, namely, the longitudinal cruise speed and the occupied lane in the highway, respectively, see Fig.\ \ref{fig:AutPlayers} for an illustration. This motivates us to model the dynamics of each noncooperative vehicle over a certain horizon as a \gls{MLD} system \cite{bemporad1999control}.

\begin{figure}[!h]
	\centering
	\includegraphics[width=.6\columnwidth]{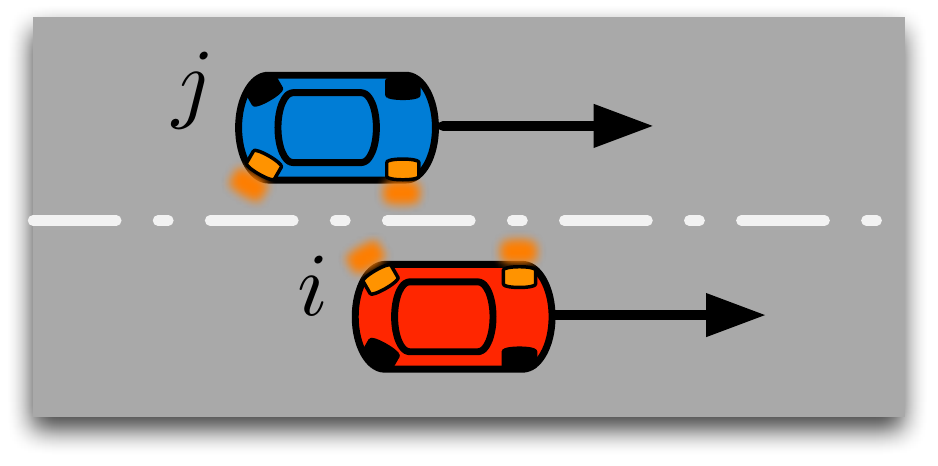}
	\caption{Two vehicles traveling side by side along two consecutive lanes. Vehicle $j$ (left lane) has the right indicator on, while $i$ (right lane) the left one.}
	\label{fig:AutDriv_FreeSpace_Lane}
\end{figure}

\begin{figure}[!h]
	\centering
	\includegraphics[width=0.95\columnwidth]{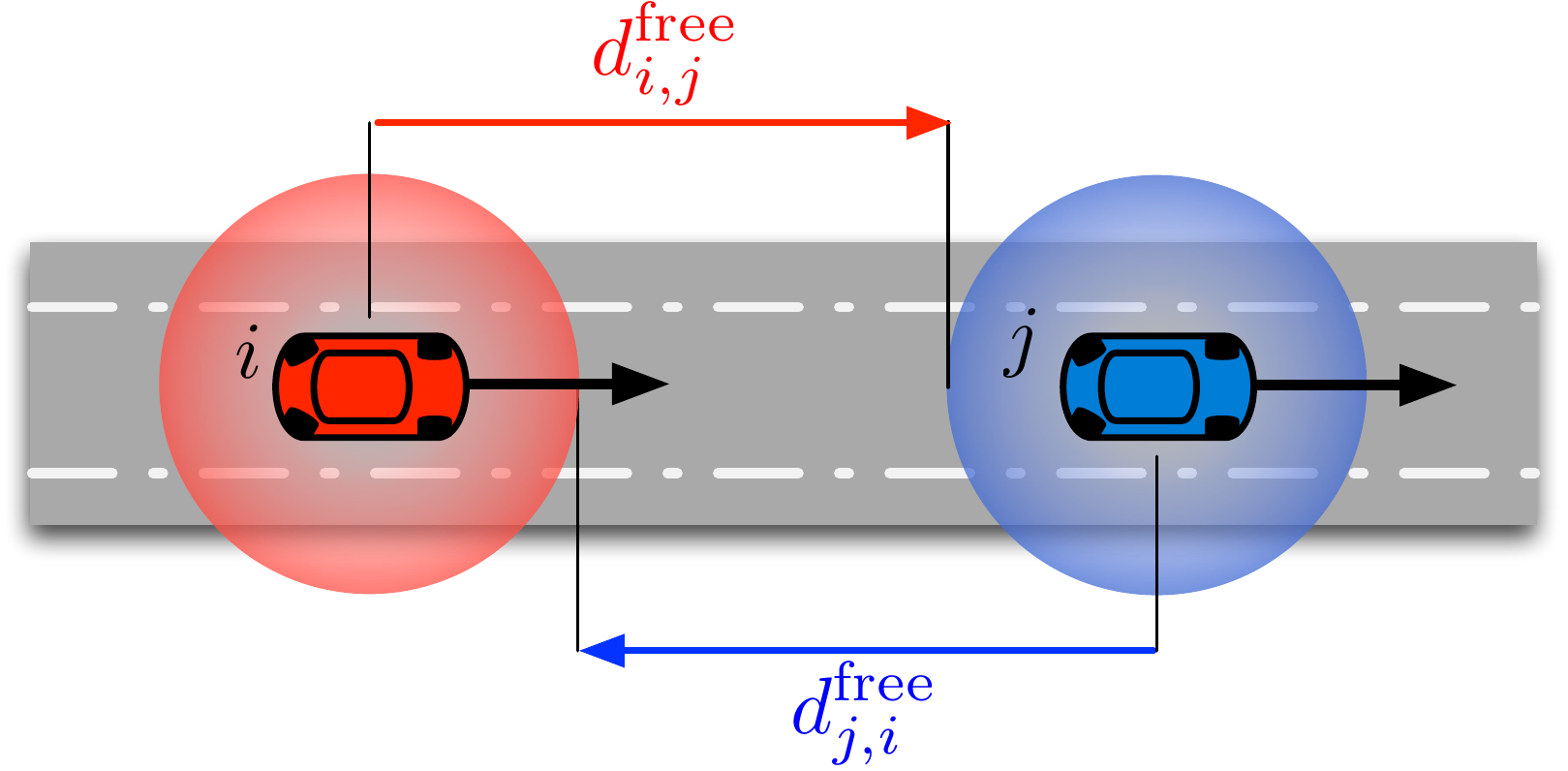}
	\caption{Two vehicles traveling on the same lane and competing for the longitudinal free space.}
	\label{fig:AutDriv_FreeSpace}
\end{figure}

The presence of multiple noncooperative agents with mixed-integer decision variables and safety constraints complicates enormously the solution of the inter-dependent decision-making problems, as conflicts naturally arise \cite{lygeros1998verified}. For instance, one conflict arises when two vehicles aim at swapping their lanes by simultaneously activating their direction indicators, see Fig.~\ref{fig:AutDriv_FreeSpace_Lane}. Conflicts arise even on a single lane, when a fast vehicle approaches a slower one, hence the two ``compete'' for the free space, see Fig.~\ref{fig:AutDriv_FreeSpace}.
 

Our approach to resolve conflicts and potential collisions is by introducing some ``AD rules'' (\S \ref{sec:MICons}). For simplicity, we assume that each vehicle is aware of the planning of its neighboring vehicles, e.g. by exchanging information and/or estimating the motion of the neighboring vehicles (\S\ref{sec:AD_MLD}). Finally, we formalize the multi-vehicle AD problem as a generalized mixed-integer potential game (\S\ref{sec:AD_MIPG}), where an equilibrium solution, obtained via a Gauss-Southwell algorithm (\S\ref{sec:Sim}), corresponds to a sequence of mixed-integer decisions for the vehicles that are individually optimal, given the decisions of the other vehicles and the imposed AD rules.

\section{Highway traffic as a system of mixed-logical-dynamical system}\label{sec:AD_MLD}

Let $\mc{I} \coloneqq \{1,\ldots,N\}$ be the set of vehicles driving on a highway with lane set $\mc{L} \coloneqq \{1,\ldots, L\}$. 
For any pair of vehicles $(i,j)$, let $d_{i,j} \in \R$ denote the inter-vehicle distance between $i$ and $j$. Throughout the paper, we refer to $i$ as a generic vehicle in $\mc{I}$ and to $j$ as a vehicle in the neighborhood of vehicle $i$, i.e., $\mc{N}_i \coloneqq \{j \in \mc{I}\, \mid \, \lvert d_{i,j} \rvert \leq \bar{d} \}$, where $\bar{d} > 0$ denotes a predefined interaction distance.
As shown in Fig.~\ref{fig:AutPlayers}, we assume that each vehicle $i$ controls its longitudinal (cruise) speed $v_i \in \mc{V}_i \subset \R$ and selects the traveling lane $z_i \in \mc{L}$. Over a prediction horizon $T \geq 1$, each vehicle $i$ has decision variables  $\bs{v}_i \coloneqq \left[v_i(1); \ldots; v_i(T)\right] \in \mc{V}_i^{T}$ and $\bs{z}_i \coloneqq \left[z_i(1); \ldots; z_i(T)\right] \in \mc{L}^{T}$. 
We assume that vehicle $i$ seeks for a sequence of hybrid decisions that trade off the tracking of a desired speed profile $\bs{v}_i^\text{d} \in \mc{V}_i^{T}$, while driving along a desired lane $\bs{z}_i^\text{d} \in \mc{L}^T$. 
Therefore, we can preliminary formulate the \gls{MPC} motion planning as a \gls{MILP}:
\begin{equation}\label{eq:MILP_i_incomplete}
\left\{\begin{aligned}
&\underset{\bs{v}_i, \bs{z}_i}{\textrm{min}} & & |\bs{v}_i - \bs{v}_i^\text{d}| + r_i |\bs{z}_i - \bs{z}_i^\text{d}|\\
&\hspace{.15cm}\textrm{s.t.} & & v_i(t+1) \in \mc{V}_i(t), \; \forall t \in \{0,\ldots,T-1\}\\
&&& z_i(t+1) \in \mc{L}_i(t), \; \forall t \in \{0,\ldots,T-1\},
\end{aligned}\right.
\end{equation}
where $r_i > 0$. 
The sets $\mc{V}_i$ and $\mc{L}_i \subset \mc{L}$ shall be defined to limit the cruise speed and its variations, and the selected lane. For instance, given $\bar{v}_i$ and $ \Delta_i > 0$ as the maximum velocity and acceleration/deceleration for vehicle $i$, respectively, we can define:
\begin{align*}
\mc{V}_i(t) & \coloneqq \left[0, \, \bar{v}_i\right] \cap \left[v_i(t) - \Delta_i,\, v_i(t) + \Delta_i\right] \\
\mc{L}_i(t) & \coloneqq \mc{L} \cap \left[z_i(t) - 1,\, z_i(t) + 1\right].
\end{align*}

To model the longitudinal distance between pairs of vehicles, we adopt the Euler forward scheme as updating rule for the relative distance between the vehicle $i$ and $j \in \mc{N}_i$:
\begin{equation}\label{eq:dist_dynam}
d_{i,j}(t+1) \coloneqq d_{i,j}(t) + \tau \left(v_j(t) - v_i(t)\right),
\end{equation}
where $\tau > 0$ denotes the length of a predefined time interval. It follows that within the introduced hybrid \gls{MPC} framework, each vehicle can estimate the relative distance with respect to its neighboring vehicles by knowing their velocities.

In this paper, we do not address the issue of communication among vehicles. Our aim is instead to design a hybrid framework capable to model the \gls{AD} problem in highways. 
Specifically, we focus on the mixed-integer decision-making layer for coordination and motion planning of the vehicles. 
Therefore, in the remainder of the paper, we assume that: i) each vehicle is driven autonomously by the solution of the hybrid decision-making framework; ii) vehicles can exchange information, i.e., their decision variables, without communication delays or packet loss. 
By starting from \eqref{eq:MILP_i_incomplete}, and in the spirit of \cite{bemporad1999control}, we introduce several mixed-logical coupling constraints among vehicles that lie within a certain set, with the aim to ensure safety.

\subsection{Safety distance}
The first mixed-logical coupling constraint we introduce refers to the safety distance among vehicles traveling on the same lane. 
Directly from the common driving experience in highways, it seems reasonable to assume the safety distance, $d_i^{\text{s}} > 0$, as a function of the actual cruise speed, $d_i^{\text{s}} = d_i^{\text{s}}(v_i(t))$. For instance, compared to driving at high speed, we are induced to get closer to the vehicle ahead at low speed.
Thus, let us define the discrete variable $l_{i,j} \coloneqq z_j - z_i$,
where $z_i$ and $z_j$ are the lane selected by vehicle $i$ and $j$, respectively, and introduce the following logical implications, for all $t \in \{0,\ldots,T\}$:
\begin{equation}\label{eq:logic_cons}
\left[l_{i,j}(t) = 0\right] \wedge \left[|d_{i,j}(t)| \geq 0\right] \implies \left[|d_{i,j}(t)| \geq d_i^{\text{s}}(t) \right].
\end{equation}
The necessary conditions on the left-hand side, which must occur simultaneously, allow to select only those vehicles that, in the prediction of vehicle $i$, occupy the same lane. 
Note that the inequality $|d_{i,j}(t)| \geq 0$ allows one to cluster the vehicles in $\mc{N}_i$ as either ahead $i$ or behind $i$. Then, for all $t$ and $j$ for which both the conditions are met, the relative distance $d_{i,j}(t)$ must be greater or equal than the safety distance $d_i^{\text{s}}(t)$.

\smallskip
\begin{definition}[Longitudinal safety]\label{def:Safe_ij}
	A pair of vehicles $(i,j) \in \mc{I}^2$ is longitudinally safe over the prediction horizon $T$ if, for all $t \in \{0, \ldots, T\}$ such that $z_i(t) = z_j(t)$, $\lvert d_{i,j}(t) \rvert \geq d_i^{\text{s}}(t)$ and, furthermore, if $z_i(t+1) = z_j(t+1)$,  $d_{i,j}(t+1) \cdot d_{i,j}(t) \geq 0$.
	The system is longitudinally safe over the prediction horizon $T$ if any pair of vehicles $(i, j) \in \mc{I}^2$ is longitudinally safe.
	\hfill$\square$
\end{definition}

\subsection{Direction indicators}
Lane change maneuvers are particularly challenging to automate because each vehicle has to adapt its actions to several road users. Inspired by a common practice in a multi-lane environment,
here we introduce integer-linear constraints to characterize the direction indicators and their utilization in the lane change maneuver. Next, we will show how to exploit it to rule out potential source of collision among vehicles, as in Fig.~\ref{fig:AutDriv_FreeSpace_Lane}. 

To model the direction indicators, we introduce two binary decision variables, $a_i^{\text{r}}$ and $a_i^{\text{l}}$. Specifically, $a_i^{\text{r}} = 1$ denotes that vehicle $i$ has its right direction indicator on, hence wants to change its current lane, moving to the right; analogously, $a_i^{\text{l}} = 1$ denotes that vehicle $i$ has left direction indicator turned on. At each time $t$, the vehicles may turn only one indicator on. This translates into an exclusive OR constraint:
\begin{equation}\label{eq:DirInd_C2}
a_i^{\text{l}}(t) + a_i^{\text{r}}(t) \leq 1.
\end{equation}

Next, we impose that vehicles may perform a lane-change maneuver only after activating the suitable direction indicator. Thus, for all $t \in \{0,\ldots,T-1\}$, we define 
\begin{equation}\label{eq:lane_cons_dir}
	\bar{\mc{L}}_i(t) \coloneqq \mc{L}_i(t) \cap \left[z_i(t) - a_i^{\text{r}}(t),\, z_i(t) + a_i^{\text{l}}(t)\right].
\end{equation}

hence, the feasible-lane constraint reads as $z_i(t+1) \in \bar{\mc{L}}_i(t)$. Although the activation of the direction indicator is mandatory before a lane change, we remark that the constraints that define $\bar{\mc{L}}_i$ do not force the lane change, but instead they make it possible in the next time interval.

\section{Safety rules for multi-lane traffic}\label{sec:Rules}
\begin{figure}[t]
	\subfloat[\label{fig:Velocities_Collision}]
	{\includegraphics[width=.5\columnwidth]{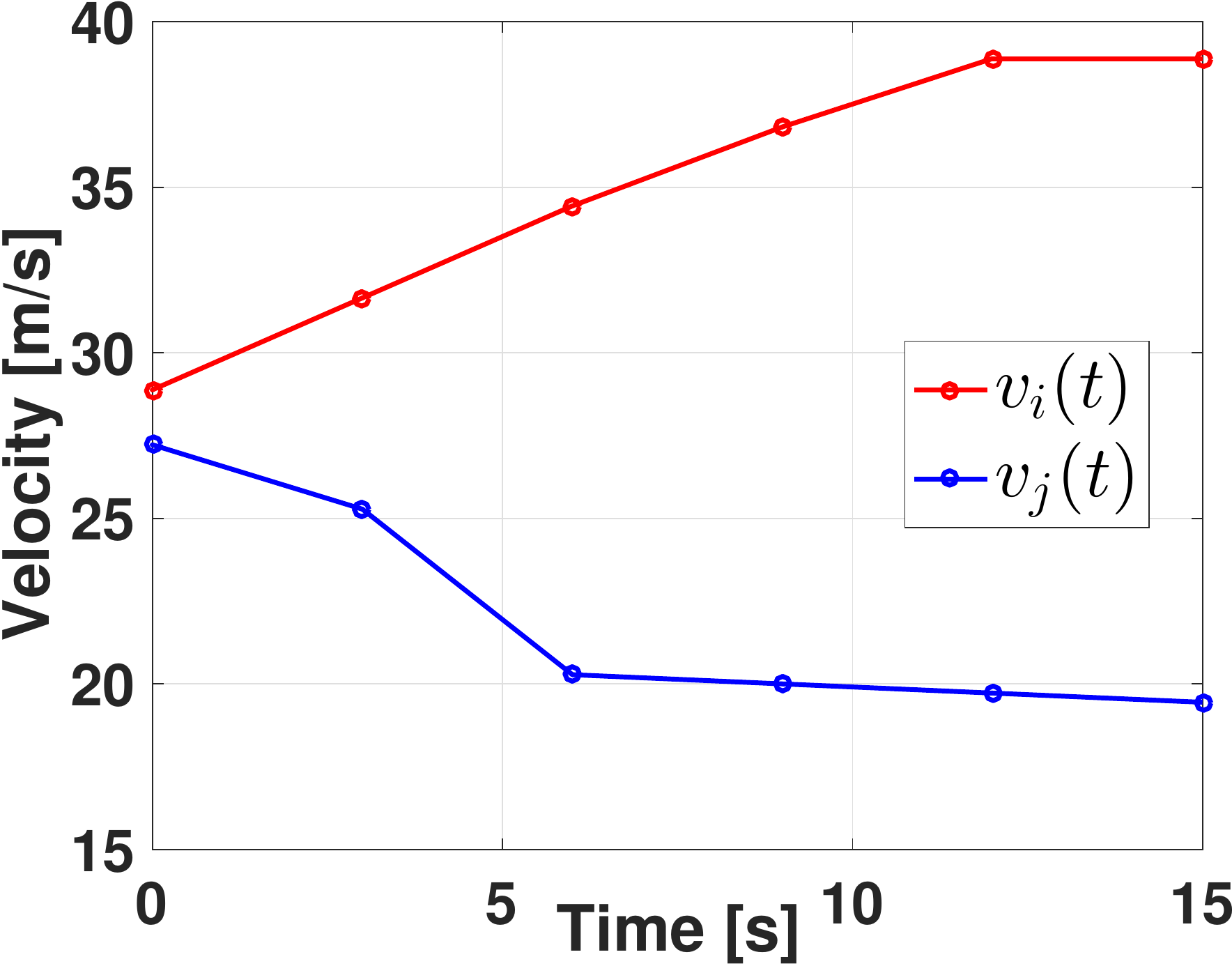}}
	\hfill
	\subfloat[\label{fig:RelDist_Collision}]
	{\includegraphics[width=.5\columnwidth]{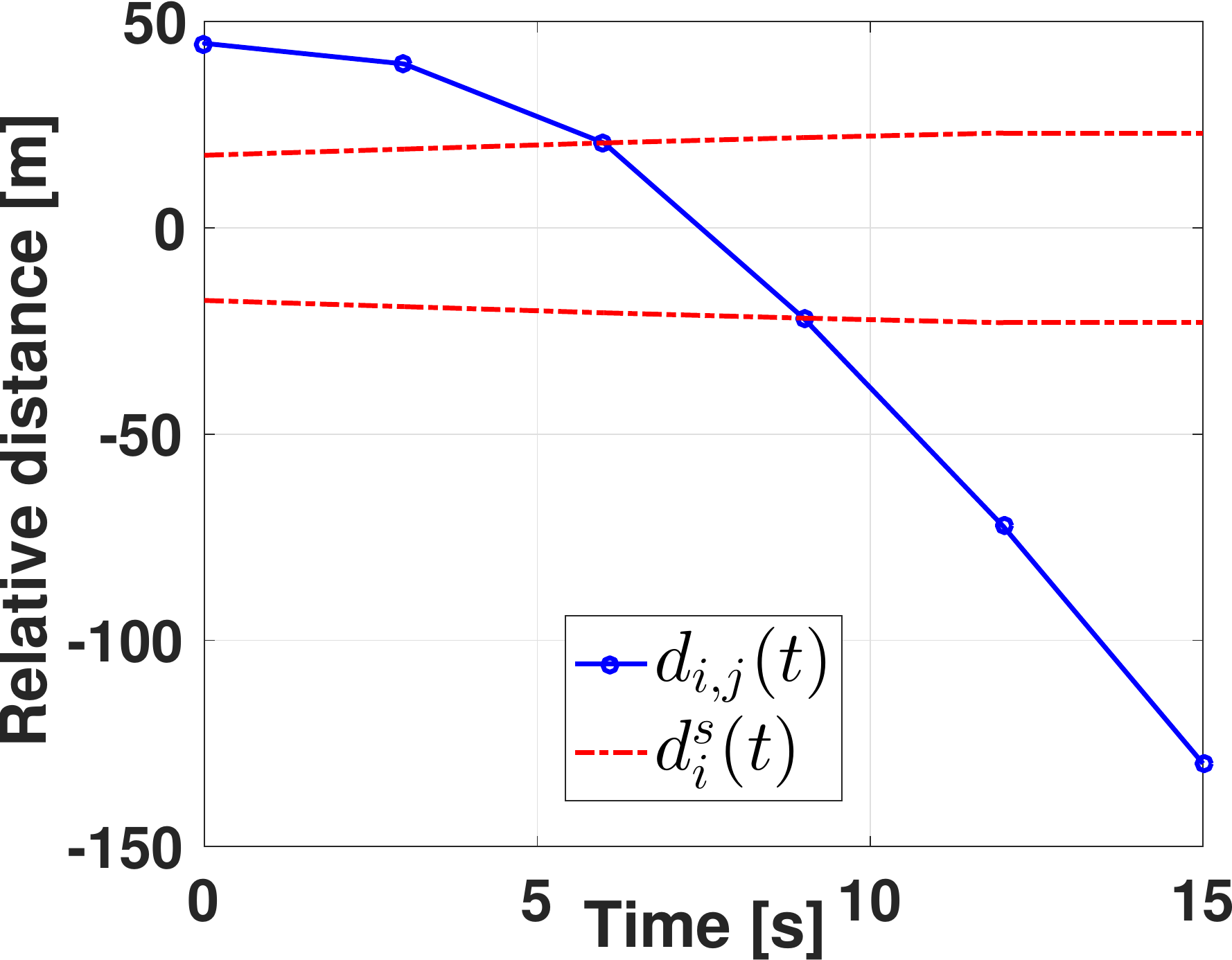}}
	\caption{Example of longitudinal collision: (a) Velocity profiles. (b) Relative and safety distances.}
	\label{fig:CollisionExample}
\end{figure}									
The logical implications introduced above allow for unsafe driving scenarios. In the following, after analyzing these scenarios, we propose some ``autonomous driving rule'' that rule out potential sources of collision.

\subsection{A free-space agreement on the lane}
Let us consider the situation in \figname~\ref{fig:AutDriv_FreeSpace}, where a feasible scenario is that vehicle $i$ accelerates, to e.g. minimize its traveling time, while vehicle $j$ reduces its speed, e.g. with the aim to minimize fuel consumption. In terms of velocity profile, an optimal strategy exists for both the vehicles (\figname~\ref{fig:Velocities_Collision}), i.e., the \glspl{MILP} problems in \eqref{eq:MILP_i_incomplete} with additional constraints are feasible. However, since the vehicles travel on the same lane, such strategies are clearly not implementable because they lead to a collision, as shown in \figname~\ref{fig:RelDist_Collision}.

To rule out this unsafe scenario, we propose an  ``agreement'' on the free space available between vehicles traveling on the same lane. Specifically, under the same necessary conditions in \eqref{eq:logic_cons}, we impose that the relative velocity between two vehicles in consecutive steps shall be limited.

In details, let us refer to \figname~\ref{fig:AutDriv_FreeSpace} and let introduce $v_{i,j} \coloneqq v_j - v_i$ as the relative velocity between vehicle $i$ and $j$. 

Then, we have, for all $j \in \mc{N}_i$ and $t \in \{1,\ldots,T\}$:
\begin{subequations}\label{eq:FreeSpace_Agg}
	\begin{align}
	\left[l_{i,j}(t) = 0\right] &\wedge \left[d_{i,j}(t) \geq 0\right] \nonumber\\ & \implies \left[v_{i,j}(t) \geq - \tfrac{1}{2}\tfrac{d_{i,j}(t) - d_i^{\text{s}}(t)}{\tau} \right],\label{eq:FreeSpaceAgg_first}\\
	\left[l_{i,j}(t) = 0\right] &\wedge \left[d_{i,j}(t) \leq 0\right] \nonumber\\ & \implies \left[v_{i,j}(t) \leq - \tfrac{1}{2} \tfrac{d_{i,j}(t) + d_i^{\text{s}}(t)}{\tau} \right].\label{eq:FreeSpaceAgg_second}
	\end{align}
\end{subequations}

Informally speaking, at each time interval, each vehicle is allowed to (selfishly) exploit only a portion (at most half) of the free longitudinal space. We emphasize that the condition $d_{i,j}(t+1) \cdot d_{i,j}(t) \geq 0$ in Definition~\ref{def:Safe_ij} would introduce nonlinear constraints. In fact, this motivates our agreement rule in \eqref{eq:FreeSpace_Agg} that introduces mixed-integer linear constraints. 

\smallskip
\begin{proposition}\label{prop:FreeSpaceAgg}
Given a pair of vehicles $(i, j) \in \mc{I}^2$, assume that $v_i(0), v_j(0)$ are feasible. The hybrid \gls{MPC} formulation in \eqref{eq:MILP_i_incomplete} with safety distance constraints \eqref{eq:logic_cons} and free-space agreement \eqref{eq:FreeSpace_Agg} guarantees the longitudinal safety.
\hfill$\square$
\end{proposition}

\smallskip
\begin{proof}
Without restriction, assume that $d_i^{\text{s}}(t) = d^{\text{s}} > 0$, for all $t \in \{1,\ldots,T\}$ and $i \in \{1,2\}$. 
The free space at step $t$ between two vehicles is $d_{i,j}(t) - d^{\text{s}} = -d_{j,i}(t) - d^{\text{s}} > 0$. Directly from \eqref{eq:FreeSpaceAgg_first}, we have $\tau v_i(t) \leq \tau v_j(t) + \tfrac{1}{2}(d_{i,j}(t) - d^{\text{s}})$. Therefore, by the definition of $d_{i,j}(t+1)$ in \eqref{eq:dist_dynam}, we obtain $d_{i,j}(t+1) - d_{i,j}(t) \geq \tfrac{1}{2}(d^{\text{s}} - d_{i,j}(t))$, which turns into $d_{i,j}(t+1) \geq \tfrac{1}{2}(d^{\text{s}} + d_{i,j}(t))$. From \eqref{eq:logic_cons}, $d_{i,j}(t) \geq d^{\text{s}}$. 
In the worst case, i.e., when $d_{i,j}(t) = d^{\text{s}}$, we obtain $d_{i,j}(t+1) \geq d^{\text{s}}$.
Now, a longitudinal collision happens if:
\begin{equation}\label{eq:lemma1_ClimbOver}
\left\{ 
\begin{aligned}
-&d_{i,j}(t+1) + d_{i,j}(t) \geq 2 d^{\text{s}}\\
&d_{i,j}(t+1) \leq 0,
\end{aligned} \right.\text{ i.e., }
\left\{ \begin{aligned}
&v_{i,j}(t) \leq - 2 \tfrac{d^{\text{s}}}{\tau}\\
&v_{i,j}(t) \leq -\tfrac{d_{i,j}(t)}{\tau}.
\end{aligned} \right.
\end{equation}
By \eqref{eq:FreeSpaceAgg_first}, we have that:
\begin{equation}\label{eq:lemma1_FreeSpace}
2 d_{i,j}(t+1) - d_{i,j}(t) \geq d^{\text{s}}\implies v_{i,j}(t) \geq \tfrac{1}{2}\tfrac{d^{\text{s}} - d_{i,j}(t)}{\tau}.
\end{equation}
Finally, \eqref{eq:lemma1_FreeSpace} fulfill the conditions in \eqref{eq:lemma1_ClimbOver} if:
\begin{equation*}
\left\{\begin{aligned}
&\tfrac{1}{2}\tfrac{d^{\text{s}} - d_{i,j}(t)}{\tau} \leq -\tfrac{d_{i,j}(t)}{\tau}\\
&\tfrac{1}{2}\tfrac{d^{\text{s}} - d_{i,j}(t)}{\tau} \leq -2\tfrac{d^{\text{s}}}{\tau} 
\end{aligned} \right.
\implies 
\left\{ \begin{aligned}
&d_{i,j}(t) \leq -d^{\text{s}}\\
&d_{i,j}(t) \geq 5 d^{\text{s}}
\end{aligned} \right.
\end{equation*}
The latter system has no solution due to the fact that $d^{\text{s}} > 0$. This implies that the free-space agreement is sufficient to avoid that conditions in \eqref{eq:lemma1_ClimbOver} may occur.
\end{proof}

\subsection{The need for direction indicators}		
\begin{figure}[t]
	\subfloat[\label{fig:LaneSel_Collision}]
	{\includegraphics[width=.5\columnwidth]{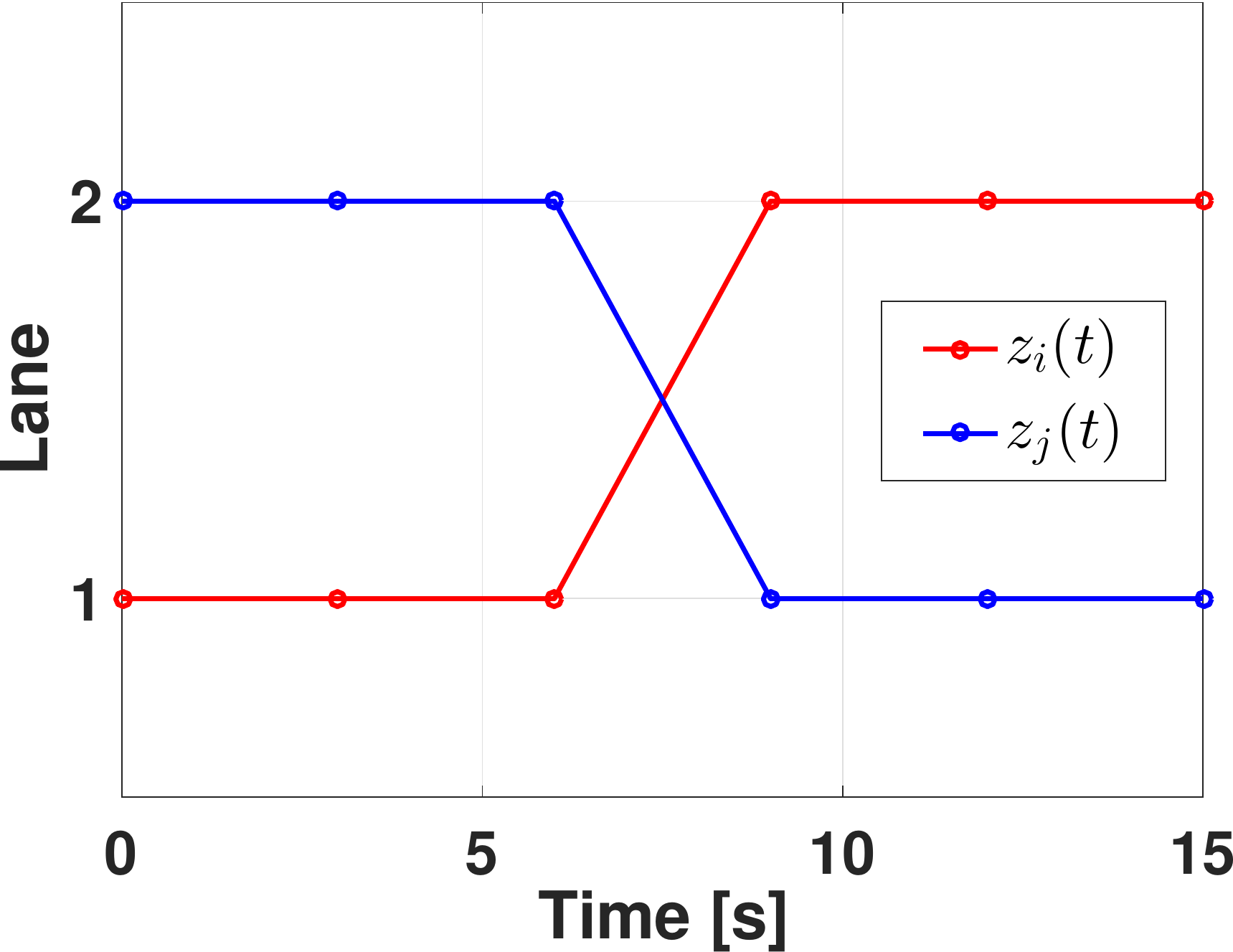}}
	\hfill
	\subfloat[\label{fig:RelDist_LaneSel_Collision}]
	{\includegraphics[width=.5\columnwidth]{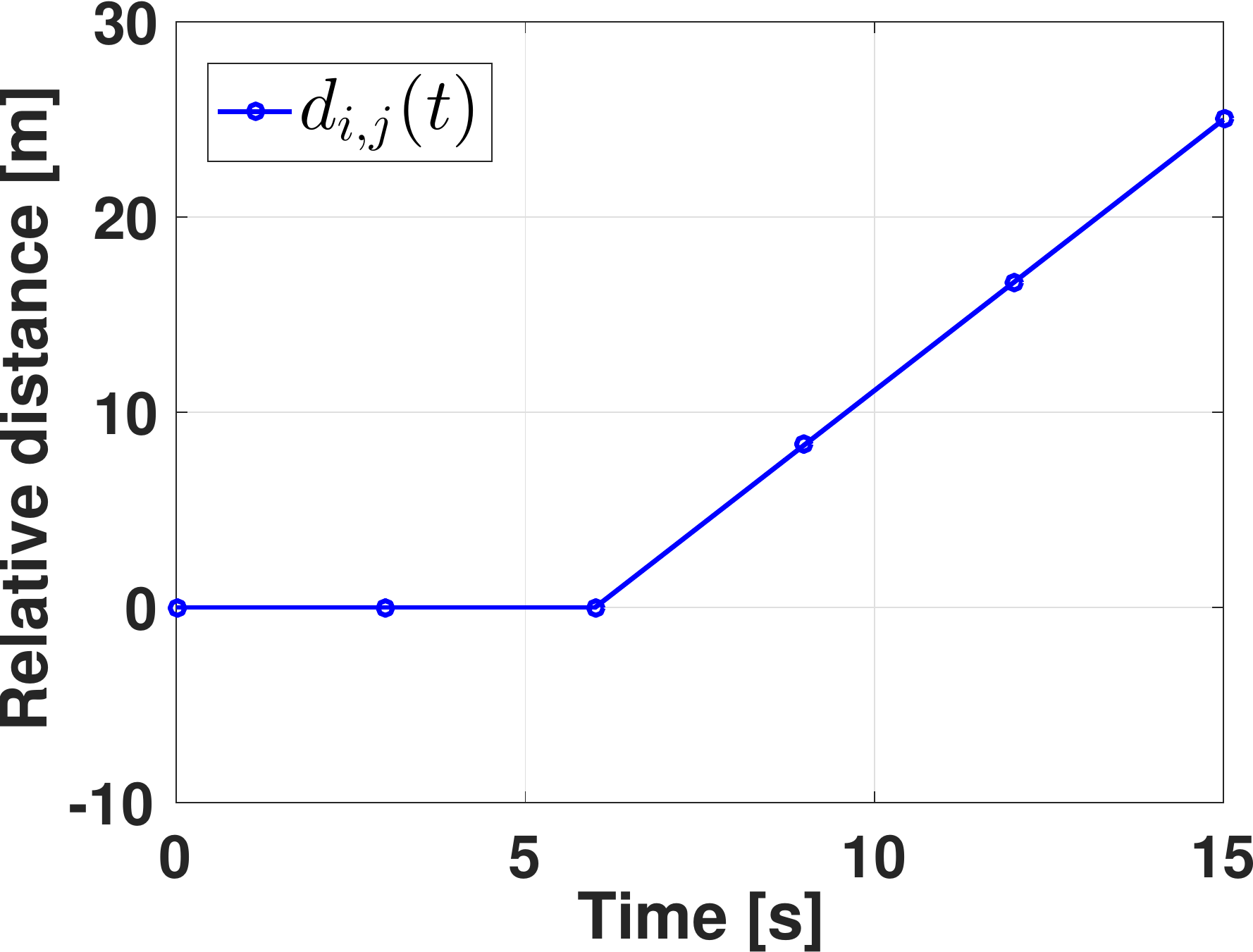}}
	\caption{Example of lateral collision between two vehicles: (a) Lane profiles. (b) Relative distance.}
	\label{fig:LateralCollision}
\end{figure}

In this subsection, we exploit the direction indicators to avoid unsafe scenarios between vehicles on consecutive lanes. Let us consider the scenario depicted in \figname~\ref{fig:AutDriv_FreeSpace_Lane}. The safety-distance logical implications in \eqref{eq:logic_cons} do not allow to change lane individually over the prediction horizon $T$, due to the small relative distance between them.  
However, the situation is different if both vehicles aim to perform the same maneuver in ``opposite directions'', swapping the lanes as showed in \figname~\ref{fig:LaneSel_Collision}. 
In this case, both predict that the destination lane will be free during the successive time intervals. 
Therefore, it is possible that by keeping their own speed unchanged, as well as relative distance, the two vehicles perform the lane change at the same time, causing a collision. However, these unsafe maneuvers are feasible for the hybrid motion planning in \eqref{eq:MILP_i_incomplete}. 

We then introduce $\hat{d} > 0$ as the inter-distance between vehicles that could lead to a lateral collision during a simultaneous change lane. We remark that $\hat{d}$ shall be chosen large enough to exclude potential conflict on consecutive lane, accordingly to the following definition.
 
\smallskip
\begin{definition}[Consecutive lane safety]\label{def:MultiSafe_ij}
	A pair of vehicles $(i,j) \in \mc{I}^2$ is safe on consecutive lanes over the prediction horizon $T$ if, for all $t \in \{0, \ldots, T\}$ such that $|d_{i,j}(t)| \leq \hat{d}$ and $|l_{i,j}(t)| = 1$, $z_i(t+1) \neq z_j(t)$ and $z_j(t+1) \neq z_i(t)$.
	The system is safe on consecutive lanes over the prediction horizon $T$ if any pair of vehicles $(i, j) \in \mc{I}^2$ is safe on consecutive lanes.
	\hfill$\square$
\end{definition}
\smallskip

To avoid lateral collisions caused by simultaneous lane changes, we propose an additional mixed-logical rule that exploits the direction indicators.
Without restriction, we refer to a scenario involving a pair of vehicles as the one illustrated in \figname~\ref{fig:AutDriv_FreeSpace_Lane}. Thus, two vehicles $(i, j)$ travel side by side on consecutive lanes if $|d_{i,j}(t)| \leq \hat{d}$ and $l_{i,j}(t) = 1$. In case of both vehicles express the will of change lane performing a swap, i.e., vehicle $i$ turns on the left indicator $a_i^{\text{l}}$, while the vehicle $j$ the right one $a_j^{\text{r}}$ at the same time $t$, then we force the vehicle traveling on a lower lane to keep it, that is $z_i(t+1) = z_i(t)$. 
Note that the proposed solution is one possible solution to resolve conflicts on consecutive lanes. Higher lanes are usually deputed for overtaking maneuvers, hence vehicles should facilitate the re-entry towards lower lanes. This motivates our proposed solution.
Thus, the logical rule reads as:
\begin{align}\label{eq:logic_implications_inDir}
&[l_{i,j}(t) = 1] \wedge [|d_{i,j}(t)| \leq \hat{d}] \wedge \{[a_i^{\text{l}}(t) = 1] \wedge [a_j^{\text{r}}(t) = 1]\}\nonumber\\
&\implies [z_i(t+1) - z_i(t) = 0]
\end{align}

\smallskip
\begin{proposition}
	Given a pair of vehicles $(i, j) \in \mc{I}^2$ and some $\hat{d} > 0$ sufficiently large, assume that $[v_i(0); z_i(0)]$ and $[v_j(0); z_j(0)]$ are feasible. The hybrid \gls{MPC} formulation in \eqref{eq:MILP_i_incomplete} with  rule \eqref{eq:logic_implications_inDir} guarantees the consecutive lane safety.
	\hfill$\square$
\end{proposition}

\smallskip
\begin{proof}
By Def.~\ref{def:MultiSafe_ij}, two vehicles might not be safe on consecutive lanes if $z_i(t+1) = z_j(t)$ and $z_j(t+1) = z_i(t)$. In view of \eqref{eq:lane_cons_dir}, this is possible only if $|l_{i,j}(t)| = 1$ and each vehicle turns on the proper direction indicator, i.e., $a_i^{\text{l}}(t) = 1$ and $a_j^{\text{r}}(t) = 1$ (or $a_i^{\text{r}}(t) = 1$ and $a_j^{\text{l}}(t) = 1$). If $|d_{i,j}(t)| > \hat{d}$ the vehicles may swap the lanes without collision; otherwise \eqref{eq:logic_implications_inDir} forces the vehicle driving on a lower lane to keep it.
\end{proof}

\section{From logical implications to mixed-integer linear constraints}\label{sec:MICons}
In this section we show how to translate the logical implications in \eqref{eq:logic_cons}, \eqref{eq:FreeSpace_Agg}, \eqref{eq:logic_implications_inDir} into mixed-integer linear constraints suitable for \eqref{eq:MILP_i_incomplete}. By referring to the vehicle $i$, we introduce the constraints to be designed for each neighboring vehicle $j \in \mc{N}_i$ and for each time $t \in \{1,\ldots, T\}$.
\subsection{Preliminaries}
Let us consider the safety distance constraints in \eqref{eq:logic_cons}. We introduce two further logical implications and related binary variables, $\alpha, \beta \in \{0,1\}$, which allow to discriminate only such vehicles that effectively travel along the same lane ($l_{i,j} = 0$) of the $i$-th one, either ahead ($\beta = 1$) or behind it ($\beta = 0$):
\begin{subequations}\label{eq:aux_logic_impl}
	\begin{align}
	&\left[\alpha_{i,j}(t) = 1\right] \iff \left[l_{i,j}(t) \leq 0\right] \wedge \left[l_{i,j}(t) \geq 0\right],\label{eq:aux_logic_impl_first}\\
	&\left[\beta_{i,j}(t) = 1\right] \iff \left[d_{i,j}(t) \geq 0\right].\label{eq:aux_logic_impl_second}
	\end{align}
\end{subequations}

Hence, equations \eqref{eq:logic_cons}, \eqref{eq:FreeSpace_Agg} can be rewritten as nonlinear inequalities:
\begin{subequations}\label{eq:LinIneq_1}
\begin{align}
\alpha_{i,j}(t) &\left[\beta_{i,j}(t) \left(d^{\text{s}}_i(t) - d_{i,j}(t)\right) \right.\nonumber\\ 
&\left. + \left(1 - \beta_{i,j}(t)\right) \left(d^{\text{s}}_i(t) + d_{i,j}(t)\right) \right] \leq 0,\label{eq:LinIneq_1_1}\\
\alpha_{i,j}(t) &\left[ -\beta_{i,j}(t) \left(2 \tau v_{i,j}(t) + d_{i,j}(t) - d^{\text{s}}_i(t)\right) \right.\nonumber\\
&\left.+ \left(1 - \beta_{i,j}(t)\right) \left(2 \tau v_{i,j}(t) + d_{i,j}(t) + d^{\text{s}}_i(t)\right)\right] \leq 0.\label{eq:LinIneq_1_2}
\end{align}
\end{subequations}

In a similar way, it is possible to handle \eqref{eq:logic_implications_inDir} with some extra binary variables, $\gamma$, $\delta$ and $\zeta$. Thus, we then introduce the associated logical implications:
\begin{subequations}\label{eq:aux_logic_impl_2}
	\begin{align}
	&\left[\gamma_{i,j}(t) = 1\right] \iff \left[l_{i,j}(t) \leq 1\right] \wedge \left[l_{i,j}(t) \geq 1\right],\label{eq:aux_logic_impl_2_first}\\
	&\left[\delta_{i,j}(t) = 1\right] \iff [a_i^{\text{l}}(t) = 1] \wedge [a_j^{\text{r}}(t) = 1],\label{eq:aux_logic_impl_2_second}\\
	& \left[\zeta_{i,j}(t) = 1\right] \iff [d_{i,j}(t) \leq \hat{d}] \vee [d_{i,j}(t) \geq -\hat{d}],\label{eq:aux_logic_impl_2_third}
	\end{align}
\end{subequations}
which allow to equivalently reformulate \eqref{eq:logic_implications_inDir} as:
\begin{equation}\label{eq:LinIneq_2}
\gamma_{i,j}(t) \delta_{i,j}(t) \zeta_{i,j}(t) \left(z_i(t+1) - z_i(t)\right) = 0.
\end{equation}

\subsection{The mixed-integer linear constraints}
For the sake of clarity, we define several patterns of inequalities that allow to handle all the constraints. Given a linear function $f:\R \rightarrow \R$, let define $M \coloneqq \max_{x \in \mc{X}} f(x)$, $m \coloneqq \min_{x \in \mc{X}} f(x)$ with $\mc{X}$ compact set. Then, with $c \in \R$ and $\delta \in \{0,1\}$, a first system $\mc{S}_{\geq}$ of mixed-integer inequalities correspond to $[\delta = 1] \iff [f(x) \geq c]$, i.e.,
\begin{equation*}
\mc{S}_{\geq}(\delta, f(x), c) \coloneqq \left\{
\begin{aligned}
&(c - m)\delta \leq f(x) - m\\
&(M - c + \epsilon) \delta \geq f(x) - c + \epsilon,
\end{aligned}
\right.
\end{equation*}
while a second $\mc{S}_{\leq}$ to $[\delta = 1] \iff [f(x) \leq c]$:
\begin{equation*}
\mc{S}_{\leq}(\delta, f(x), c) \coloneqq \left\{
\begin{aligned}
&(M - c) \delta \leq M - f(x)\\
&(c + \epsilon - m) \delta \geq \epsilon + c - f(x).
\end{aligned}
\right.
\end{equation*}

Here $\epsilon > 0$ is a small tolerance beyond which the constraint is regarded as violated. As an example, let consider the right-hand side in \eqref{eq:aux_logic_impl_first}: introducing $\eta$, $\theta \in \{0,1\}$, $\left[\eta_{i,j}(t) = 1\right] \iff \left[l_{i,j}(t) \leq 0\right]$ translates into $\mc{S}_{\leq}(\eta_{i,j}(t), l_{i,j}(t), 0)$, while $\left[\theta_{i,j}(t) = 1\right] \iff \left[l_{i,j}(t) \geq 0\right]$ into $\mc{S}_{\geq}(\theta_{i,j}(t), l_{i,j}(t), 0)$. 
Moreover, we define the next two blocks of inequalities, involving only binary variables, which allow to solve propositions with logical AND:
\begin{equation*}
\mc{S}_{\wedge}(\delta, \sigma, \gamma) \coloneqq \left\{
\begin{aligned}
- &\sigma + \delta \leq 0\\
- &\gamma + \delta \leq 0\\
&\sigma + \gamma - \delta \leq 1,
\end{aligned}
\right.
\end{equation*}
and with logical OR:
\begin{equation*}
\mc{S}_{\vee}(\delta, \sigma, \gamma) \coloneqq \left\{
\begin{aligned}
&\sigma - \delta \leq 0\\
&\gamma - \delta \leq 0\\
- &\sigma - \gamma + \delta \leq 0.
\end{aligned}
\right.
\end{equation*}

Specifically, $[\delta = 1] \iff [\sigma = 1] \wedge [\gamma = 1]$ is equivalent to the integer inequalities $\mc{S}_{\wedge}(\delta, \sigma, \gamma)$, while $[\delta = 1] \iff [\sigma = 1] \vee [\gamma = 1]$ into $\mc{S}_{\vee}(\delta, \sigma, \gamma)$. Referring again to \eqref{eq:aux_logic_impl_first}, $[\alpha_{i,j} (t) = 1] \iff [\eta_{i,j}(t) = 1] \wedge [\theta_{i,j}(t) = 1]$ corresponds to $\mc{S}_{\wedge}(\alpha_{i,j}(t), \eta_{i,j}(t), \theta_{i,j}(t))$. Finally, \eqref{eq:aux_logic_impl_first} coincides with the system of mixed-integer inequalities given by:
\begin{equation}\label{eq:IneqSys_1}
\eqref{eq:aux_logic_impl_first} \Longrightarrow \left\{
\begin{aligned}
&\mc{S}_{\leq}(\eta_{i,j}(t), l_{i,j}(t), 0),\\
&\mc{S}_{\geq}(\theta_{i,j}(t), l_{i,j}(t), 0),\\
&\mc{S}_{\wedge}(\alpha_{i,j}(t), \eta_{i,j}(t), \theta_{i,j}(t)).
\end{aligned} \right.
\end{equation}

Thus, it follows that:
\begin{align}
\eqref{eq:aux_logic_impl_second} &\Longrightarrow \mc{S}_{\geq}(\beta_{i,j}(t), d_{i,j}(t), 0).\\
\eqref{eq:aux_logic_impl_2_first} &\Longrightarrow \left\{
\begin{aligned}
&\mc{S}_{\leq}(\kappa_{i,j}(t), l_{i,j}(t), 1),\\
&\mc{S}_{\geq}(\lambda_{i,j}(t), l_{i,j}(t), 1),\\
&\mc{S}_{\wedge}(\gamma_{i,j}(t), \kappa_{i,j}(t), \lambda_{i,j}(t)).
\end{aligned}\right.\\
\eqref{eq:aux_logic_impl_2_second} &\Longrightarrow \mc{S}_{\wedge}(\delta_{i,j}(t), a^{\text{l}}_{i}(t), a^{\text{r}}_{j}(t)).\\
\eqref{eq:aux_logic_impl_2_third} &\Longrightarrow \left\{
\begin{aligned}
&\mc{S}_{\leq}(\mu_{i,j}(t), d_{i,j}(t), \hat{d}),\\
&\mc{S}_{\geq}(\nu_{i,j}(t), d_{i,j}(t), -\hat{d}),\\
&\mc{S}_{\vee}(\zeta_{i,j}(t), \mu_{i,j}(t), \nu_{i,j}(t)).
\end{aligned}\right.
\end{align}

Next, we follow the procedure in \cite{bemporad1999control} to recast the inequalities in \eqref{eq:LinIneq_1} and \eqref{eq:LinIneq_2} into a mixed-integer linear formulation by means of additional auxiliary variables (both real and binary, \cite{williams2013model}). Specifically, starting from \eqref{eq:LinIneq_1}, we define $\xi_{i,j} \coloneqq \alpha_{i,j} \beta_{i,j}$, which satisfies the system of inequalities
\begin{equation}
\mc{S}_{\wedge}(\xi_{i,j}(t), \alpha_{i,j}(t), \beta_{i,j}(t)).
\end{equation} 

By referring to \eqref{eq:LinIneq_1_1}, we also define the real auxiliary variables $f_{i,j} \coloneqq \xi_{i,j} d_{i,j}$, $g_{i,j} \coloneqq \alpha_{i,j} d^{\text{s}}_i$ and $h_{i,j} \coloneqq \alpha_{i,j} d_{i,j}$ that shall satisfy the pattern of linear inequalities given by:
\begin{equation*}
\mc{S}_{\Rightarrow}(g, f(x), \delta) \coloneqq \left\{
\begin{aligned}
& m \delta \leq g \leq M \delta\\
& - M(1-\delta) \leq	g - f(x) \leq - m(1 - \delta)\\
\end{aligned}\right.
\end{equation*}

The latter is equivalent to: $[\delta = 0] \implies [g = 0]$, while $[\delta = 1] \implies [g = f(x)]$. Hence, for each real auxiliary variable previously introduced, we have the systems:
\begin{align}
&\mc{S}_{\Rightarrow}(f_{i,j}(t), d_{i,j}(t), \xi_{i,j}(t)),\\
&\mc{S}_{\Rightarrow}(g_{i,j}(t), d^{\text{s}}_i(t), \alpha_{i,j}(t)),\\
&\mc{S}_{\Rightarrow}(h_{i,j}(t), d_{i,j}(t), \alpha_{i,j}(t)).
\end{align}

Thus, the nonlinear inequalities in \eqref{eq:LinIneq_1_1} becomes:
\begin{equation}\label{eq:SecDist_Lin}
-2 f_{i,j}(t) + g_{i,j}(t) + h_{i,j}(t) \leq 0.
\end{equation}

Now, let us consider \eqref{eq:LinIneq_1_2}. We define two real auxiliary variables, $k_{i,j} = \xi_{i,j} v_i$ and $m_{i,j} = \alpha_{i,j} v_i$, that satisfy:
\begin{align}
&\mc{S}_{\Rightarrow}(k_{i,j}(t), v_i(t), \xi_{i,j}(t)),\\
&\mc{S}_{\Rightarrow}(m_{i,j}(t), v_i(t), \alpha_{i,j}(t)).
\end{align}

Hence, \eqref{eq:LinIneq_1_2} is rewritten with linear formulation as:
\begin{equation}\label{eq:FreeSpace_Lin}
\begin{aligned}
2 \tau (2 k_{i,j}(t) &- m_{i,j}(t)) - 2 f_{i,j}(t) + g_{i,j}(t) + h_{i,j}(t)\\&+ 2 \tau (\alpha_{i,j}(t) - 2 \xi_{i,j}(t)) v_j(t)  \leq 0.
\end{aligned}
\end{equation}

Finally, we proceed with the same procedure as for \eqref{eq:LinIneq_2} by introducing two auxiliary binary variables, $\phi_{i,j} \coloneqq \gamma_{i,j} \delta_{i,j}$ and $\psi_{i,j} \coloneqq \phi_{i,j} \zeta_{i,j}$, that satisfy the systems
\begin{align}
&\mc{S}_{\wedge}(\phi_{i,j}(t), \gamma_{i,j}(t), \delta_{i,j}(t)),\\
&\mc{S}_{\wedge}(\psi_{i,j}(t), \zeta_{i,j}(t), \phi_{i,j}(t)),
\end{align}

and two discrete variables, $p_{i,j} \coloneqq \psi_{i,j} z_i(t)$ and $s_{i,j} \coloneqq \psi_{i,j} z_i(t+1)$, so that we obtain:
\begin{equation}\label{eq:LaneInd_Lin}
\left\{
\begin{aligned}
-&s_{i,j}(t) + p_{i,j}(t) \leq 0\\
& s_{i,j}(t) - p_{i,j}(t) \leq 0.
\end{aligned} \right.
\end{equation}

Then, the variables $s_{i,j}$ and $u_{i,j}$ satisfy the inequalities
\begin{align}
&\mc{S}_{\Rightarrow}(p_{i,j}(t), z_i(t), \psi_{i,j}(t)),\\
&\mc{S}_{\Rightarrow}(s_{i,j}(t), z_i(t+1), \psi_{i,j}(t)).\label{eq:IneqSys_2}
\end{align}

\subsection{The final mixed-integer linear model}
In the previous subsection, for each vehicle $i$, we have introduced 21 auxiliary variables, both continuous and discrete, and 71 mixed-integer linear constraints. By rearranging all the inequalities, we propose the final \gls{MILP} for each vehicle: 
\begin{equation}\label{eq:MILP_i_complete}
\left\{\begin{aligned}
&\underset{q_i, \bs{v}_i, \ldots,\bs{s}_i}{\textrm{min}} & & q_i\\
&\hspace{.5cm}\textrm{s.t.} & & -q_i \bs{1} \leq \bs{v}_i - \bs{v}_i^{\text{d}} \leq q_i \bs{1}\\
&&& -q_i \bs{1} \leq r_i (\bs{z}_i - \bs{z}_i^{\text{d}}) \leq q_i \bs{1}\\
&&& v_i(t+1) \in \mc{V}_i(t), \; \forall t \in \{0,\ldots,T-1\}\\
&&& z_i(t+1) \in \bar{\mc{L}}_i(t), \; \forall t \in \{0,\ldots,T-1\}\\
&&& \eqref{eq:DirInd_C2}, \; \forall t \in \{0,\ldots,T-1\}\\
&&& \eqref{eq:IneqSys_1}-\eqref{eq:IneqSys_2}, \; \forall j \in \mc{N}_i, \, \forall t \in \{1,\ldots,T\}.
\end{aligned}\right.
\end{equation}

The total number of mixed-integer linear constraints for player $i$ is $c_i \coloneqq T (67 |\mc{N}_i| + 9)$, while for the whole neighborhood is $c \coloneqq (\sum_{j \in \mc{N}_i}^{} c_j) + c_i$.
Note that the coupling constraints in $\eqref{eq:IneqSys_1}-\eqref{eq:IneqSys_2}$ contain the strategies of the neighbors as affine, given terms. Thus, by defining $\bs{x}_i \coloneqq [q_i; \bs{v}_i; \ldots ; \bs{s}_i] \in \R^{n_i}$, where $n_i \coloneqq 1 + T (21 |\mc{N}_i| + 4)$, and $\bs{x} \in \R^{n}$, $n \coloneqq (\sum_{j \in \mc{N}_i}^{} n_j) + n_i$, as the vector of all the decision variables in the neighborhood $\mc{N}_i$:
\begin{equation}\label{eq:MILP_i_compact}
\underset{\bs{x}_i}{\textrm{min}} \  w_i^\top \bs{x}_i \ \textrm{ s.t. } \ A \bs{x} \leq b
\end{equation}

for suitable $w_i \in \R^{n_i}$, $A \in \R^{c \times n}$, $b \in \R^{c}$ vectors and matrices of suitable structure.

\section{Automated Driving as a Generalized Mixed-Integer Potential Game}\label{sec:AD_MIPG}
Within our hybrid framework, selfish road users can be driven by a set of mutually influencing mixed-integer strategies obtained by solving \eqref{eq:MILP_i_compact} for all $i \in \mc{I}$. Thus, we aim at designing suitable sequences of decision variables that control each vehicle towards its own goal, without compromising the overall safety. To achieve such a trade-off, we propose to formalize the \gls{AD} coordination problem as a generalized mixed-integer potential game.

We preliminary define the feasible set of each player, namely $\mc{X}_{i}(\bs{x}_{-i}) \coloneqq \{ \bs{x}_i \in \R^{n_i} \mid A [\bs{x}_i ; \, \bs{x}_{-i}] \leq b \}$, and $\mc{X} \coloneqq \{\bs{x} \in \R^n \mid A \bs{x} \leq b\}$.
Furthermore, by noticing that each $J_{i}(\bs{x}_{i}) \coloneqq w_i^\top \bs{x}_i$ depends only on the local variable $\bs{x}_i$,  we introduce the function $P(\bs{x}) \coloneqq \sum_{i \in \mc{I}}^{} J_{i}(\bs{x}_{i})$. By \cite{facchinei2011decomposition}, $P$ is an exact potential function for the proposed \gls{AD} game because it satisfies, for all $i \in \mc{I}$, for all $\bs{x}_{-i}$, and for all  $\bs{x}_i$, $\bs{y}_i \in \mc{X}_i(\bs{x}_{-i})$,
\begin{equation*}
P(\bs{x}_i,\bs{x}_{-i}) - P(\bs{y}_i,\bs{x}_{-i}) = J_i(\bs{x}_i) - J_i(\bs{y}_i).
\end{equation*}

Let us now introduce the mixed-integer best response mapping for player (i.e., vehicle) $i$, given the strategies of its neighbors $\bs{x}_{-i}$:
\begin{equation*}
	\bs{x}^{\star}_{i}(\bs{x}_{-i}) \ \in \ 
	\textrm{arg} \underset{\bs{x}_i}{\textrm{min}} \; J_i(\bs{x}_i) \quad \textrm{s.t.} \quad (\bs{x}_i, \bs{x}_{-i}) \in \mc{X}
\end{equation*}

\smallskip
\begin{definition}[$\varepsilon$-Mixed-Integer Nash Equilibrium]
	Let $\varepsilon > 0$. $\bar{\bs{x}} \in \mc{X}$ is an \gls{MINE} of the game if, for all $i \in \mc{I}$,
	\begin{equation*}
		J_i(\bar{\bs{x}}_i) - J_i(\bs{x}^{*}_i) \leq \varepsilon,
	\end{equation*}
	where $\bs{x}^{*}_i \in \bs{x}^{\star}_{i}(\bar{\bs{x}}_{-i})$.	
	\hfill$\square$
\end{definition}
\smallskip

Any $\varepsilon$-global minimizer of the potential function $P$, i.e., any $\bar{\bs{x}} \in \mc{X}$ such that $P(\bar{\bs{x}}) \leq P(\bs{x}) + \varepsilon$ for all $\bs{x} \in \mc{X}$, is an \gls{MINE} of the generalized mixed-integer potential game \cite[Th. 2]{sagratella2017algorithms}. The converse does not hold in general.

It follows that an \gls{MINE} is a vector of (individually) optimal strategies that allows to safely coordinate a set of noncooperative vehicles driving on a highway.
\section{Game resolution and numerical simulations via Gauss-Southwell Algorithm}\label{sec:Sim}
In this section, we show numerical results obtained by solving the mixed-integer potential game associated with the \gls{AD} problem. 
To compute an \gls{MINE}, we adopt the best-response-based \gls{GS} method, described next. At each algorithmic step $k$, an arbitrary $i = i_k \in \mc{I}$ is chosen; then, vehicle $i$ updates its decision variable as follows:
\begin{equation}
	\bs{x}_i(k+1) = \left\{
		\begin{aligned}
			& \bs{x}_i^{*}(k) && \text{if } J_{i}(\bs{x}_{i}(k)) - J_{i}(\bs{x}_i^{*}(k)) \geq \varepsilon \\
			&\bs{x}_i(k) && \text{otherwise},
		\end{aligned}\right.
\end{equation}
with $\bs{x}_i^{*}(k) \in \bs{x}^{\star}_{i}(\bs{x}_{-i}(k))$, while $\bs{x}_{-i}(k+1) = \bs{x}_{-i}(k)$. The iteration goes on until $J_{i}(\bs{x}_{i}(k)) - J_{i}(\bs{x}_i^{*}(k)) \geq \varepsilon$ holds for all $i = i_k \in \mc{I}$. 
Under suitable conditions on the sequence of $i_k \in \mc{I}$, this algorithm converges to an \gls{MINE} in a finite number of steps \cite[Th. 4]{sagratella2017algorithms}.

\begin{figure*}[!t]
	\centering
	\subfloat[\label{fig:frame_3_1}]
	{\includegraphics[width=.8\columnwidth]{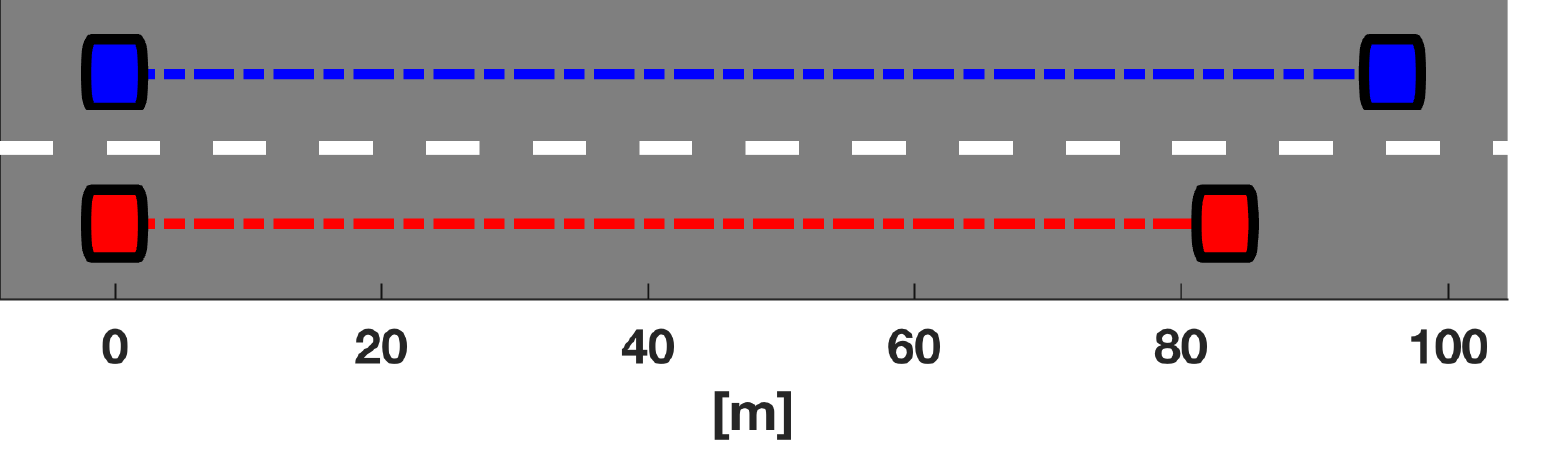}}\\
	\subfloat[\label{fig:frame_3_2}]
	{\includegraphics[width=.8\columnwidth]{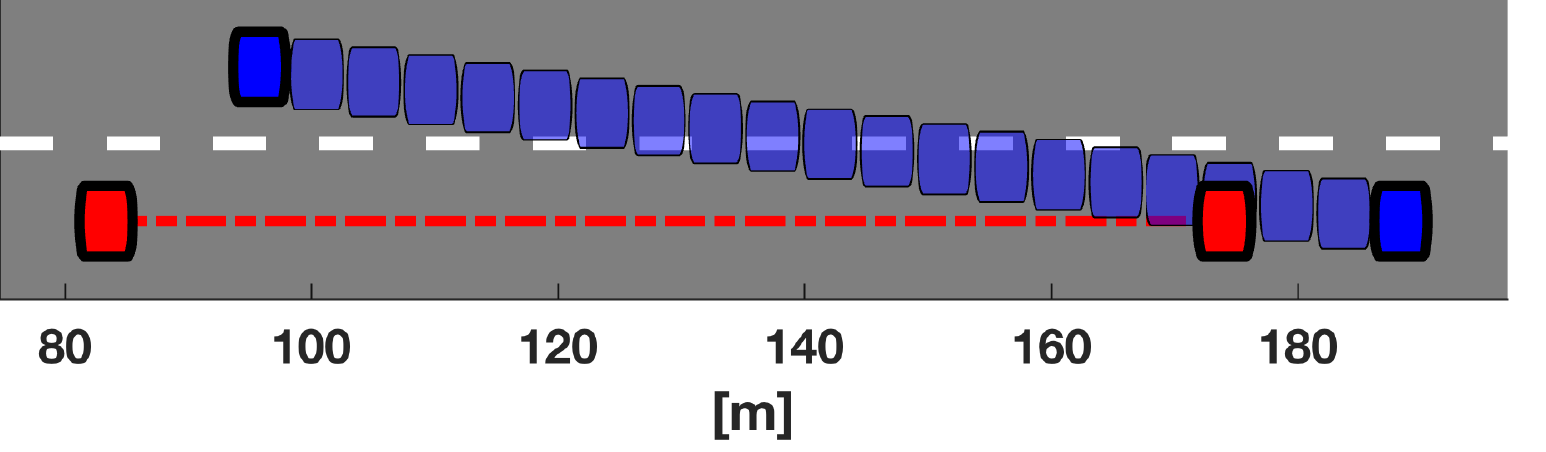}}
	\hfill
	\subfloat[\label{fig:frame_3_3}]
	{\includegraphics[width=.8\columnwidth]{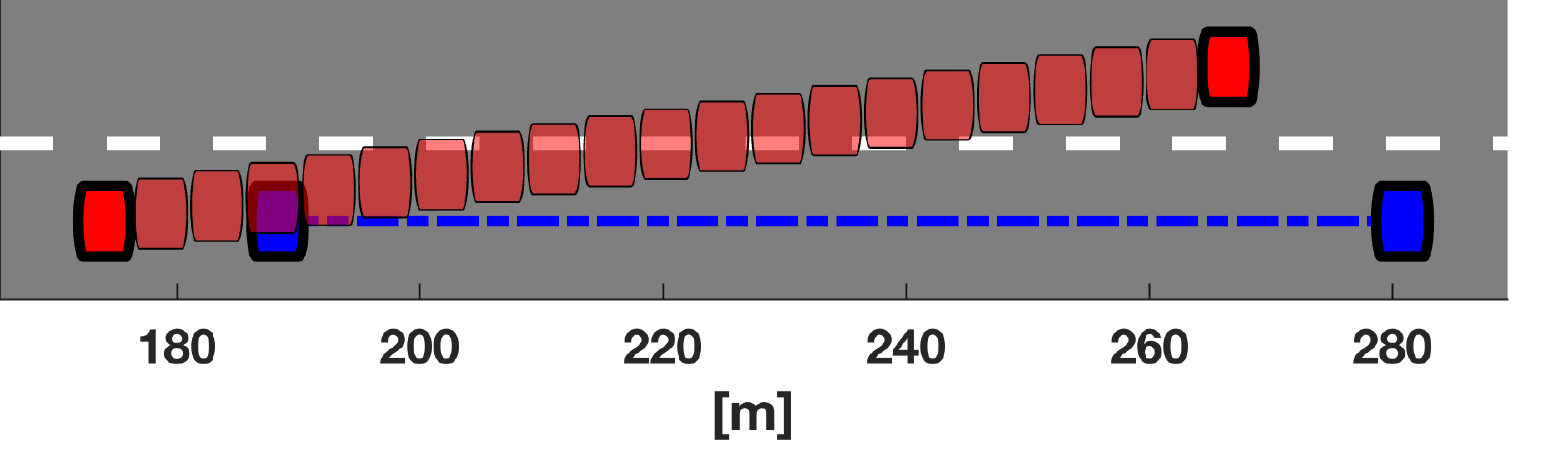}}
	\caption{Lateral conflict resolution with logic rule in \eqref{eq:logic_implications_inDir} ($T = 3$, $\tau = 3$ [s]). (a) 0--3 [s]. (b) 3.1--6 [s]. (c) 6.1--9 [s].}
	\label{fig:sim_3}
\end{figure*}

\begin{figure*}[!t]
	\centering
	\subfloat[\label{fig:sim_1_SL_vel}]
	{\includegraphics[width=.6\columnwidth]{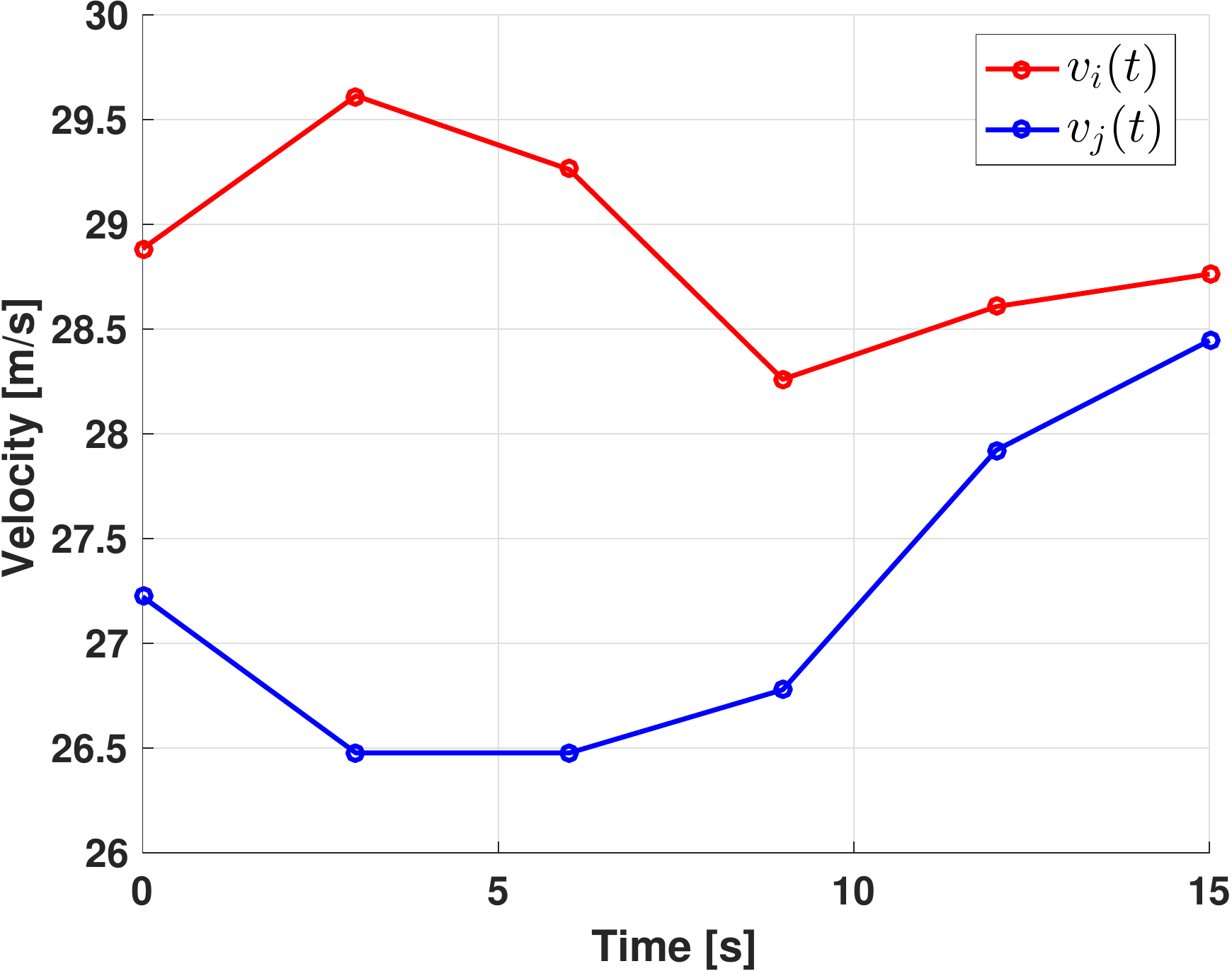}}
	\hspace{2cm}
	\subfloat[\label{fig:sim_1_SL_reldist}]
	{\includegraphics[width=.6\columnwidth]{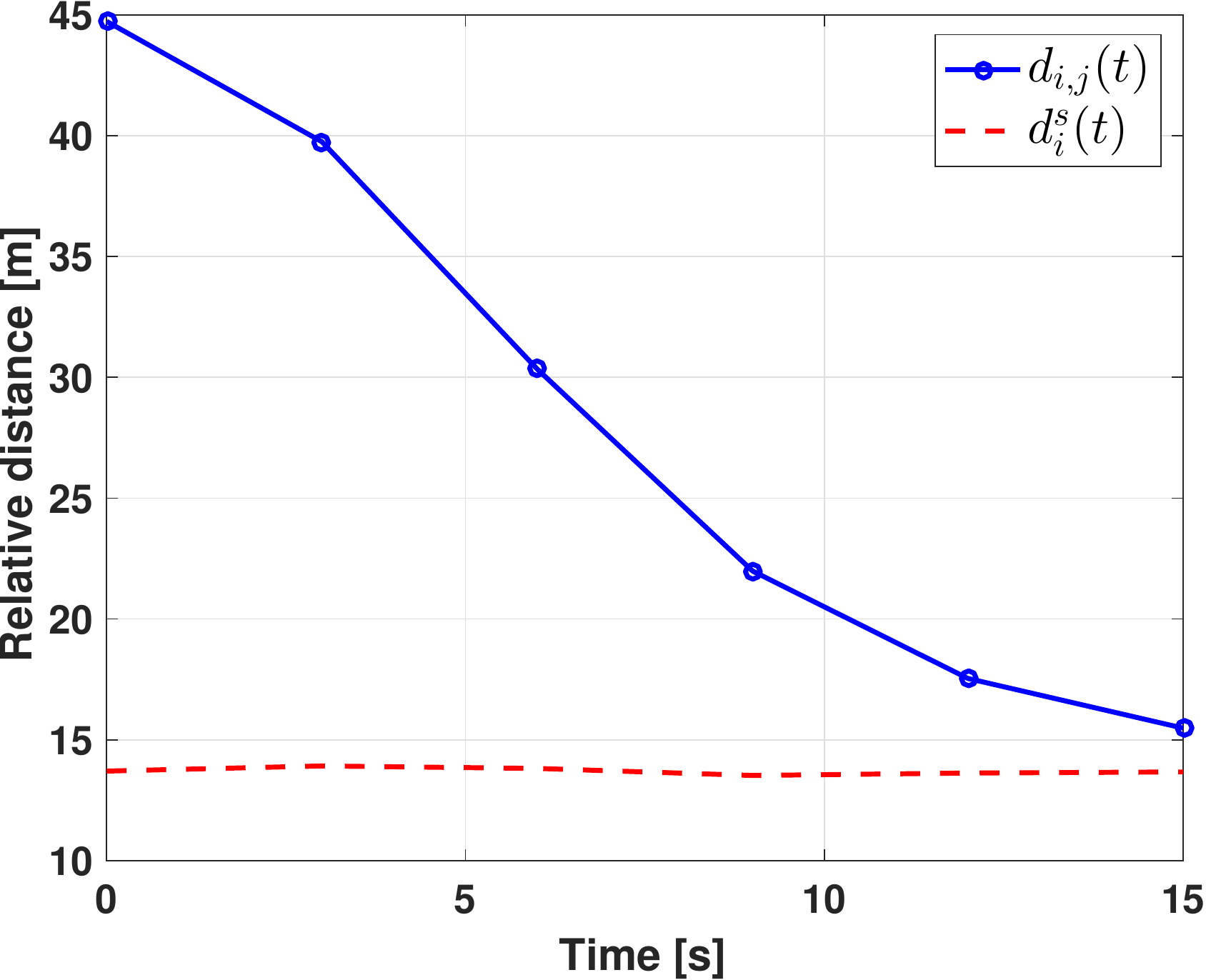}}
	\caption{Longitudinal conflict resolution with free-space agreement in \eqref{eq:FreeSpace_Agg}. (a) Velocity profiles. (b) Relative distance.}
	\label{fig:SL}
\end{figure*}

\begin{figure*}[!t]
	\centering
	\subfloat[\label{fig:frame_1}]
	{\includegraphics[width=.8\columnwidth]{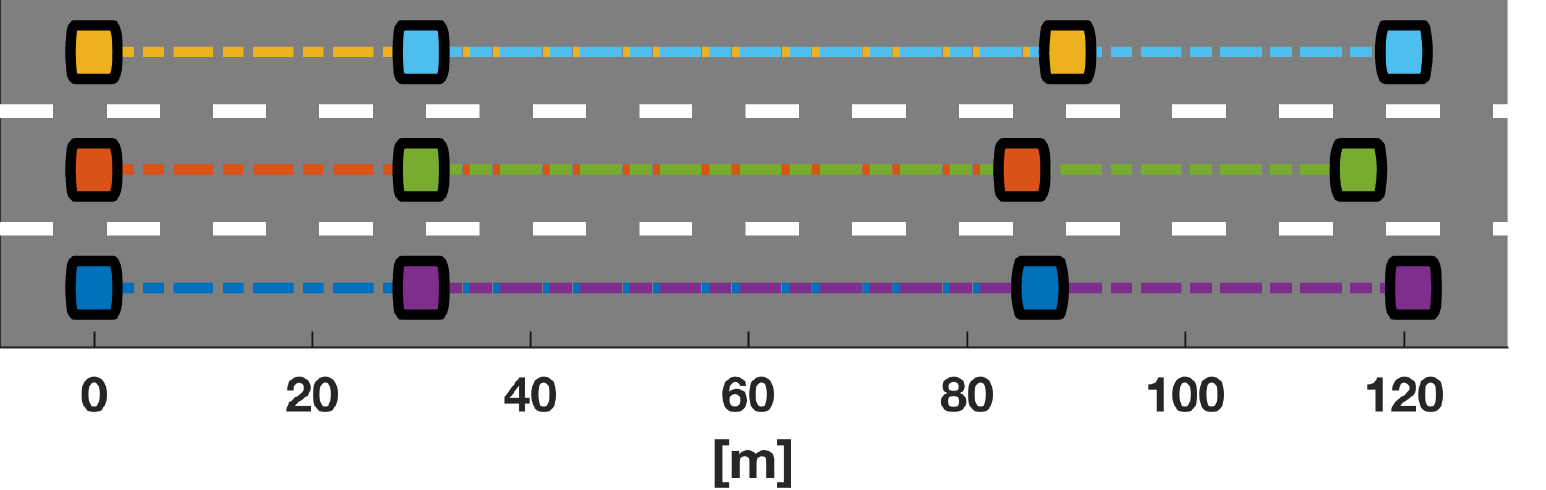}}
	\hfill
	\subfloat[\label{fig:frame_2}]
	{\includegraphics[width=.8\columnwidth]{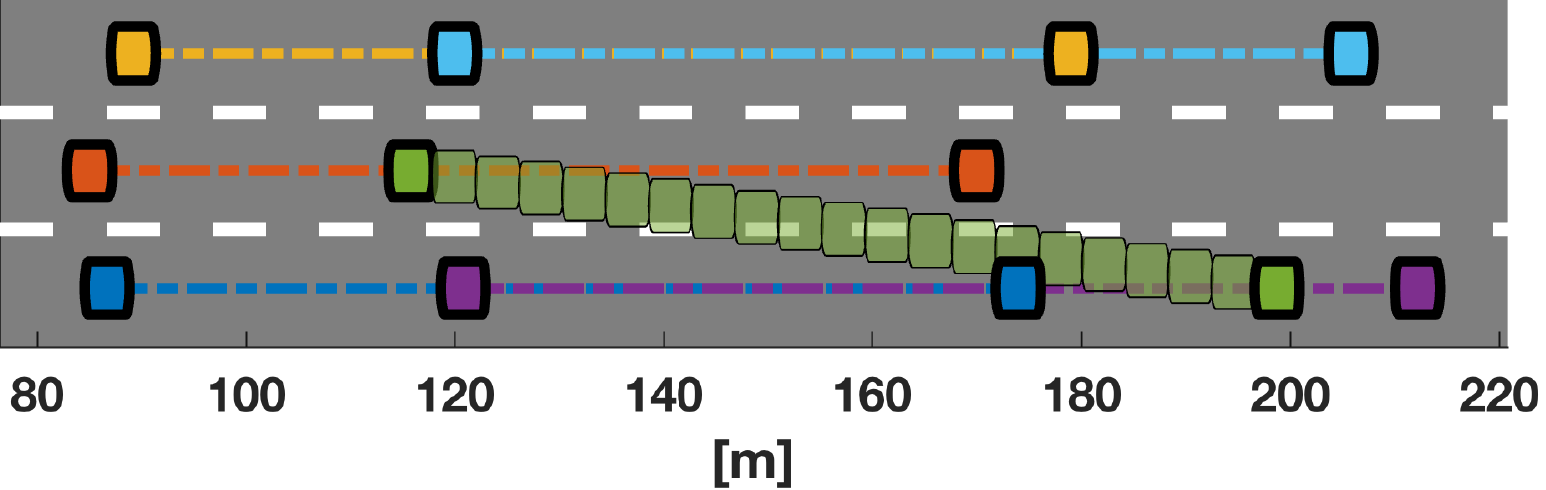}}\\
	\subfloat[\label{fig:frame_3}]
	{\includegraphics[width=.8\columnwidth]{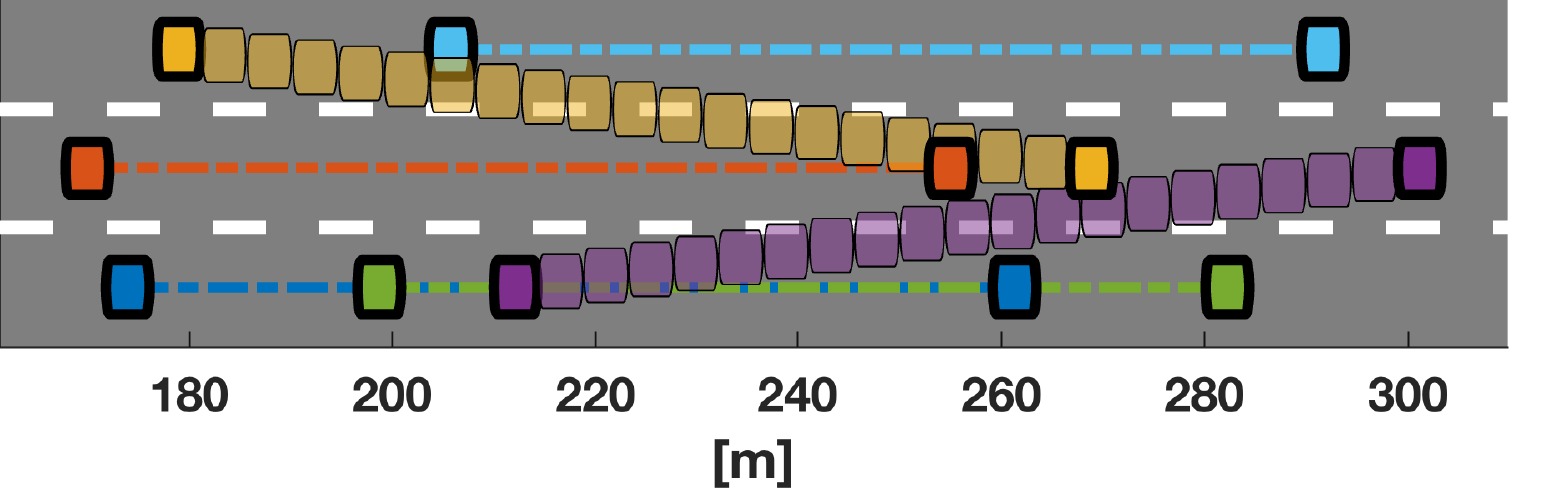}}
	\hfill
	\subfloat[\label{fig:frame_4}]
	{\includegraphics[width=.8\columnwidth]{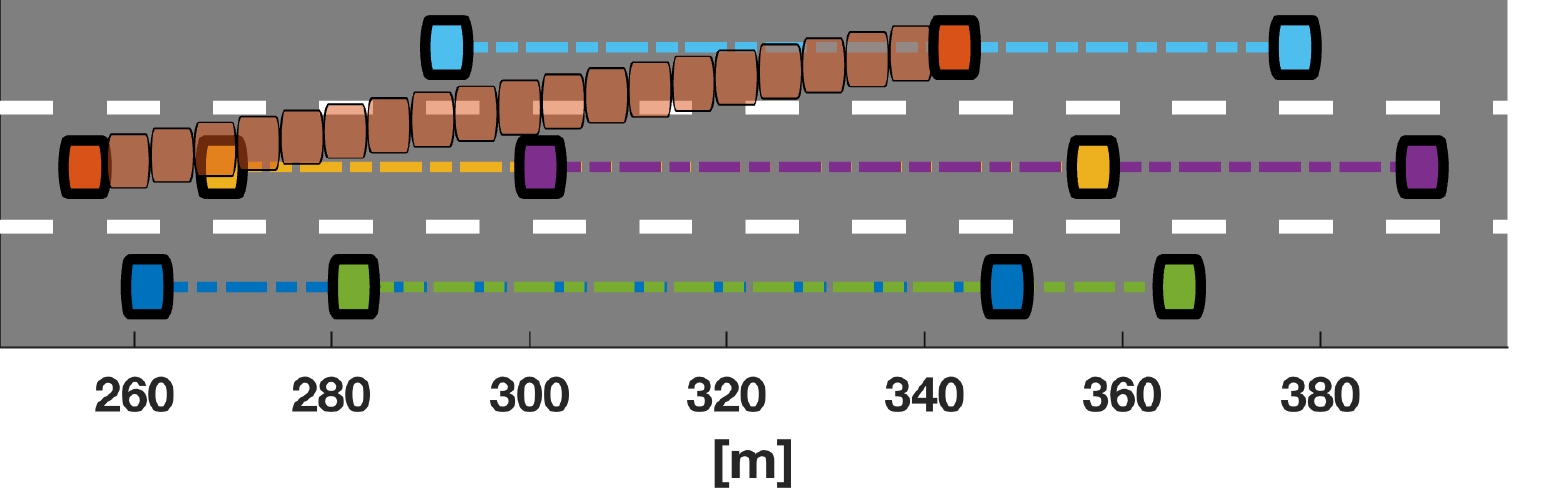}}
	\caption{Multi-lane traffic simulation with six vehicles disposed side by side on three lanes. Each vehicle achieves a randomly chosen target lane and tracks a reference speed ($T = 4$, $\tau = 3$ [s]). (a) 0--3 [s]. (b) 3.1--6 [s]. (c) 6.1--9 [s]. (d) 9.1--12 [s].}
	\label{fig:sim_10}
\end{figure*}

In Figures~\ref{fig:sim_3}--\ref{fig:SL} the logical rules allow to prevent the potential conflicts illustrated in Figures.~\ref{fig:AutDriv_FreeSpace_Lane}--\ref{fig:AutDriv_FreeSpace}, while Fig.~\ref{fig:sim_10} shows an example of multi-lane traffic simulation. Finally, Tab.~\ref{tab:data_summary} shows the average computational times for solving a single optimization problem over the full horizon $T$.

\section{Conclusion and outlook}
A hybrid decision-making framework can model the multi-lane, multi-vehicle automated driving problem in highways and, if augmented with simple driving rules, can ensure a safe use of the road space-time, despite the presence of selfish vehicles. The decision-making problem can be in fact modeled as a generalized mixed-integer potential game, which can be solved iteratively via a Gauss-Southwell best-response algorithm.
Future research will focus on closed-loop control for the generalized mixed-integer potential game that arises in multi-vehicle automated driving.


\begin{table*}[]
	\caption{Computational-time summary.}\label{tab:data_summary}
	\centering
	\begin{tabular}{>{\centering\arraybackslash}p{0.075\linewidth}>{\centering\arraybackslash}p{0.075\linewidth}>{\centering\arraybackslash}p{0.1\linewidth}>{\centering\arraybackslash}p{0.22\linewidth}}
	\toprule
	\textbf{Fig.~\#}		&	\textbf{\# vehicles}	&				\textbf{\# iterations}		&	 \textbf{Avg. computational time [ms] (to solve \eqref{eq:MILP_i_complete}) }\\
	\midrule
	\ref{fig:sim_3}		&	2	    &			4		&		15.1			\\
	\midrule
	\ref{fig:SL}		&	2	    &			4		&		14.6			\\
	\midrule
	\ref{fig:sim_10}	&	6       &			18		&		56.3		\\
	\midrule
	-	        &	9       &					27		&		69.2		\\
	\bottomrule
	\end{tabular}
\end{table*}

\bibliographystyle{IEEEtran}
\bibliography{18_CDC_AutDriv}

\end{document}